\documentclass[10pt]{article}

\usepackage{amsmath}
\usepackage{amssymb}
\usepackage{times}
\usepackage{graphicx}
\usepackage{hhline}
\usepackage{verbatim}
\usepackage{float}
\usepackage{color}
\newsavebox\dotbox
\newlength{\dotheight}
\DeclareMathOperator*{\bigcdot}{\boldsymbol{\cdot}\rule[-\dp\dotbox]{0pt}{\dotheight}}

\usepackage{amssymb}
\usepackage[dvipsnames]{xcolor}
\newcommand{\Keywords}[1]{\par\noindent
{\small{\em Keywords\/}: #1}}

\newtheorem{thm}{Theorem}
\newtheorem{lem}[thm]{Lemma}
\newtheorem{prop}[thm]{Proposition}

\newtheorem{rem}[thm]{Remark}

\newenvironment{proof}[1][Proof]{\textbf{#1.} }{\ \rule{0.5em}{0.5em}}

\def\bt{\hbox{$\bf \cdot$}}

\let\eps=\varepsilon
\let\a=\alpha
\let\b=\beta
\let\D=\Delta

\let\s=\sigma

\let\L=\Lambda

\let\ls=\leqslant
\let\gs=\geqslant

\let\cal=\mathcal

\def\dr{\hbox{\rm d}}

\def\pr{\hbox{\bf P}}
\def\E{\hbox{\bf E}}

\newcommand{\tend}[1]{\longrightarrow}
\newcommand{\tendsb}{\xrightarrow{a.s.}}

\newcommand{\tends}[1]{\xrightarrow[#1]{}}
\newcommand{\tendsd}{\xrightarrow{\ d\ }}
\newcommand{\abs}[1]{\lvert #1\rvert}

\title{On comparison of the estimators of the Hurst index and the diffusion coefficient of the fractional Gompertz diffusion process}
\date{}
\author{K. Kubilius *$^1$,  D. Melichov $^2$}

\begin{document}

\maketitle

\vskip-0.5cm

\centerline{\small $^1$ Vilnius University, Institute of Mathematics
and Informatics, Akademijos 4,} \centerline{\small LT-08663 Vilnius,
Lithuania}

\centerline{\small $^2$ Vilnius Gediminas Technical University, Faculty of Fundamental Sciences, Saul\.etekio al. 11,} \centerline{\small LT-10223, Vilnius,
Lithuania}

\let\oldthefootnote\thefootnote
\renewcommand{\thefootnote}{\fnsymbol{footnote}}
\footnotetext[1]{This research was funded by a grant (No.
MIP-048/2014) from the Research Council of Lithuania.}

\abstract{We study some estimators of the Hurst index and {the diffusion coefficient} of the fractional Gompertz diffusion process and prove that they are strongly consistent and {most} of them are asymptotically normal. Moreover, we compare the asymptotic behavior of these estimators with the aid of computer simulations.

\bigskip
 }

\section{Introduction}

Many applications make use of processes that are described by stochastic
differential equations (SDEs). Recently, much attention has been
paid to SDEs driven by the fractional Brownian motion (fBm) and to
the problems of statistical estimation of model parameters.
Statistical aspects of the models driven by the fBm have been studied in many
articles. Especially much attention has been paid to the estimation of the parameters of drift. We focus on estimators of the Hurst index and {the diffusion coefficient}. Recently some new estimators of the Hurst index and of {the diffusion coefficient} have been proposed {(see \cite{bs}, \cite{BLL}, \cite{ks}, \cite{ksm})}. This paper aims to compare them using discrete observations of the sample paths of the solution of the SDE.\\

\noindent As the test process we will consider the fractional Gompertz diffusion process (fGd)
\begin{equation}\label{e:diflygt}
X_t=\int_0^t(\alpha X_s-\beta X_s\ln X_s)\,ds+\sigma \int_0^t X_s \,dB^H_s,\qquad X_0=x_0>0,\quad 0\ls t\ls T,
\end{equation}
where  $\a$, $\b\neq 0$, and $\s$ are real parameters and $B^H$ is a fBm with the Hurst index $H\in(1/2,1)$. Almost all sample paths of $B^H$ have bounded
$p$-variation for each $p>1/H$ on $[0,T]$ for every $T>0$.  The second integral in (\ref{e:diflygt}) is the pathwise Riemann-Stieltjes integral with respect to the process having finite $p$-variation. \\

\noindent { {The reasons we have chosen fGd as the test process are as follows. Firstly, it is a} non-linear process. To the equation (\ref{e:diflygt}) {it is possible to} apply a pathwise approach and use a chain rule for the composition of a smooth function and a function of bounded $p$-variation with $1<p<2$. This approach allows to easily obtain the unique explicit solution of the equation (\ref{e:diflygt})  for $H\in(1/2,1)$ in the class of processes{, almost all sample paths of which} have bounded $p$-variation with $1<p<2$. {Secondly, the} structure of {the} increments of fGd allows us to apply {a wider class} of estimators without {imposing} additional restrictions {on the} process. The normalization of quadratic variation by {the} square of {the} process value at {a} fixed point allow{s} us to {derive the asymptotic}  normality of these estimators. {The application} of this approach allows to consider similar statistics for the equations with time-dependent coefficients.
Moreover, {in case of the} standard Brownian motion, i.e. for $H=1/2$, this process {plays an important role} in the modeling of population growth. }\\

\noindent Dung \cite{dung} proved that a class of fractional geometric mean reversion processes expressed by a fractional SDE of the form
\[
X_t=\int_0^t(\alpha_s X_s-\beta_s X_s\ln X_s)\,ds+ \int_0^t\sigma_s X_s \,d W^H_s,\qquad X_0=x_0>0,\quad 0\ls t\ls T,
\]
where $W^H$ is a fractional Brownian motion of the Liouville form, has a unique solution. It follows from his results that, if the coefficients in the equation above are constant, its solution will be of the form
\begin{equation}\label{e:sprendinys}
X_t=\exp\bigg\{e^{-\b t}\ln x_0 +\a\int_0^t e^{-\b(t-s)}ds  +\s \int_0^t e^{-\b(t-s)}\, dW^H_s\bigg\}.
\end{equation}
{In {the} Appendix} it will be shown that the equation (\ref{e:diflygt}) has the solution of the same form even without the assumption required by Dung.\\

\noindent In case of the fractional Ornstein-Uhlenbeck process and the geometric Brownian motion a comparison of various estimators of the Hurst index was presented in \cite{kubmel}. The behavior of the estimators based on quadratic variations was compared with that of some of the other known estimators. {It should be noted that these estimators are not asymptotically normal. Moreover, {only one of the estimators considered in the aforementioned paper is included in the comparison presented in this article.}}\\

\noindent A reader interested in the existence of the solution of the Gompertz diffusion process with respect to the standard Brownian motion and the estimation of its parameters is encouraged to read \cite{skiadas}, \cite{gns}, \cite{fbpl} and the references therein.\\

\noindent The structure of the paper is as follows.
Section 2 presents the estimators considered in the rest of the paper. Section 3 contains the numerical comparison of the estimators' performance. Sections 4-6 are dedicated to proofs of strong consistency of the considered estimators in case of the fractional Gompertz diffusion process. {In {the} Appendix the existence and uniqueness of {the} solution of equation (\ref{e:diflygt}) is proved.}

\section{Estimators}

In the rest of the paper we will deal with the problem of estimating the Hurst index and {the diffusion coefficient} of the fractional Gompertz diffusion process based on discrete observations of its sample paths. The estimation of the trend parameters $\a$ and $\b$, although not included in the present paper, can be performed using the least squares method. Using the change of variable $Z_t=\ln X_t$ the equation (\ref{e:diflygt}) can be reduced to the fractional Vasicek model, to which the least squares method is then applied  (see, e.g., \cite{xzzc}).

\subsection{Hurst index estimators}

Let $\pi_n=\{\tau^{m_n}_k,\ k=0,\dots,m_n\}$, $n\ge1$, $\mathbb{N}\ni
m_n\uparrow\infty$, be a sequence of partitions of the interval $[0,T]$. If  partition $\pi_n$ is uniform  then $\tau^{m_n}_k=\frac{kT}{m_n}$ for  all $k\in\{0,\ldots,m_n\}$. If $m_n\equiv n$, we write $t_k^n$ instead of $\tau_k^{m_n}$. Let $(X_t)_{t\in[0,T]}$ be a stochastic process and
\begin{gather*}
\D^{(1)}X(\tau^{m_n}_k)=X(\tau^{m_n}_k)-X(\tau^{m_n}_{k-1}),\qquad\D^{(2)} X(\tau^{m_n}_k)=X(\tau^{m_n}_{k+1})-2X(\tau^{m_n}_k) +X(\tau^{m_n}_{k-1}),\\
\qquad k=1,\ldots,m_n-(i-1), \quad i=1,2.
\end{gather*}
Denote
\[
{V^{(i)}_{m_n,T}=\sum_{k=1}^{m_n-(i-1)} \left(\frac{\D^{(i)} X(\tau^{m_n}_k)}{X(\tau^{m_n}_{k-1})}\right)^2,} \qquad i=1,2,
\]
and
\[
W_{n,k}=\sum_{j=-k_n+1}^{k_n-1} {\big(\D^{(2)} X_{s^n_j+t^n_k}\big)^2 } =\sum_{j=-k_n+1}^{k_n-1}\left(X_{s^n_{j+1}+t^n_k}-2X_{s^n_j+t^n_k} +X_{s^n_{j-1}+t^n_k}\right)^2,
\]
where $s^n_j=\frac{jT}{m_n}$, $m_n=n k_n$,  and  $k_n=n^2$.

\begin{thm}\label{t:case1}
Assume that $X$ is a solution of the fractional Gompertz SDE and $1/2<H<1$. Then
\[
\widehat H^{(j)}_n\tendsb{} H,\qquad j=1,2,3,4,
\]
and
\begin{align*}
2\ln 2\,\sqrt{n}(\widehat H^{(1)}_n-H)&\tendsd
\mathcal{N}(0,\sigma_*^2),\qquad \sigma_*^2(H)=\frac 32 \sigma^2(H)-2\sigma_{1,2}(H),\\
2\sqrt{n}\,\ln\frac{n}{T}\,(\widehat H^{(2)}_n-H)&\tendsd
N(0;\sigma^2_H),\\
\sqrt{n}(\widehat H^{(3)}_n-H)&\tendsd \mathcal{N}\bigg(0,\sigma_\ell^2\bigg(\mathbf{r},\frac 12(\mathbf{z}/\sqrt{\mathbf{r}}\,)\bigg)\bigg)
\end{align*}
with known variances $\sigma^2_H$, $\sigma_{1,2}(H)$, $\sigma_\ell^2\big(\mathbf{r},\frac 12(\mathbf{z}/\sqrt{\mathbf{r}}\,)\big)$  defined in section \ref{s:fbm}, where
\begin{align}
\widehat H^{(1)}_n=&\frac{1}{2}-\frac{1}{2\ln2}
\ln\left(\frac{V^{(2)}_{2n,T}}{V^{(2)}_{n,T}}\right),\\
\widehat H_n^{(2)}=&\frac{1}{2}+\frac{1}{2\ln k_n}\ln\left(\frac 2n
\sum_{k=2}^n\frac{\big(\D^{(2)} X(t^n_k)
\big)^2}{W_{n,k-1}}\right),\\
\widehat H_n^{(3)}=&-\frac{1}{2}\sum_{j=1}^\ell z_j \ln \bigg(\frac{V^{(2)}_{n_j,T}}{n_j-1} \bigg),\qquad n_j=r_j n,\  r_j\in \mathbb{N},\  j=1,\ldots,\ell,\ \ell\neq 1,\\
z_i=&\frac{y_i}{\sum_{i=1}^\ell y_i^2}\quad\mbox{and}\quad y_i=\ln r_i-\frac 1\ell \sum_{i=1}^\ell \ln r_i,\nonumber\\
\widehat H_n^{(4)}=&\frac{1}{0.1468} \left(\frac{1}{n^4-2}\sum_{k=1}^{n^4-2}\frac{\left\vert \D^{(2)} X(t^{n^4}_k) + \D^{(2)} X(t^{n^4}_{k+1})\right\vert}
{\left\vert\D^{(2)} X(t^{n^4}_k)\right\vert + \left\vert\D^{(2)} X(t^{n^4}_{k+1})\right\vert} - 0.5174\right).
\end{align}

\end{thm}

\begin{rem} The estimators $\widehat H^{(i)}_n$, $i=1,2,3,4$, were considered in \cite{ksm}, \cite{km}, \cite{BLL}, and \cite{bs}. The estimator $\widehat H^{(2)}_n$ can be used to estimate the Hurst index of the generic form of the SDE with an additional restriction on the diffusion coefficient.
\end{rem}

\subsection{Diffusion coefficient estimators}

In this section, we describe four estimators of {the diffusion coefficient}. The application of the fourth is not explicitly justified, however this can be performed. It was proposed in \cite{BLL} for the fractional geometric Brownian motion. The aforementioned paper shows it to be a weakly consistent estimator of {the diffusion coefficient $\s^2$.}

\begin{thm}\label{t:case2} Assume that $X$ is a solution of the fractional Gompertz SDE, $1/2<H<1$,  and $\widehat H_n=H+o_\omega(\phi(n))$, where $o_\omega$ is defined in subsection \ref{s:gompertz}. If $\phi(n)=o\left(\frac{1}{\ln n}\right)$ then
\begin{align*}
\widehat \s^{\,2}_{1,n}=&\frac{n^{2 \widehat H_n-1}}{T^{2 \widehat
H_n}}\,V^{(1)}_{n,T}\tendsb{} \s^2,\\
\widehat \s^{\,2}_{2,n}=&\frac{n^{2 \widehat H_n-1}}{T^{2 \widehat
H_n}(4-2^{2 \widehat H_n})}\,V^{(2)}_{n,T}\tendsb{} \s^2, \\
\widehat{\s}^{\,2}_{3,n}=&\frac{\sum_{k=1}^n\big(\D^{(1)} X(t^n_k)
\big)^2}{(\tfrac Tn)^{2\widehat H_n} \sum_{k=1}^n X^2(t^n_{k-1})}\tendsb{} \s^2.
\end{align*}
If $\phi(n)=o\left(\frac{1}{\sqrt{n} \ln n}\right)$ then
\begin{align*}
&\sqrt n(\widehat \s^{\,2}_{2,n}-\s^2)\tendsd {\mathcal{N}\big(0;\s^4\sigma^2(H)\big)},\qquad\mbox{for}\quad 0<H<1,
\end{align*}
where variance $\sigma^2(H)$ are defined in section \ref{s:fbm}.
\end{thm}
For the purposes of comparison we shall also consider
\[
\widehat \s_{4,n}=\frac{\exp(\widehat B)}{(4-2^{2\widehat H_n^{(3)}})}\,,\qquad \widehat B=\frac 12\bigg(\frac 1\ell \sum_{i=1}^\ell \ln \bigg(\frac{V^{(2)}_{n_i,T}}{n_i-1} \bigg)\bigg) +\widehat H^{(3)}_n \bigg(\frac 1\ell \sum_{i=1}^\ell \ln (n_i)\bigg),
\]
where $n_i$ and $\widehat H^{(3)}_n$ are defined in Theorem \ref{t:case1}.

\begin{rem} The estimators $\widehat \s^{\,2}_{i,n}$, $i=1,2$, are similar to the estimators used in the book \cite{BLL} for the evaluation of {the diffusion coefficient} $\s$ of the solutions of linear SDE when $H$ is known.  The estimator $\widehat \s^{\,2}_{3,n}$ is used to estimate {the diffusion coefficient} of the fractional Ornstein-Uhlenbeck process when $H$ is known (see \cite{xzzc}). We have shown that this restriction can be lifted.
\end{rem}

\section{Modeling of the estimators}

The goal of this section is to describe the numerical simulations that were performed in order to compare the behavior of the estimators considered in this paper.\\

\noindent The sample paths of the fractional Brownian motion, which were further used to construct the sample paths of the fractional Gompertz diffusion process, were simulated using the Wood-Chan circulant matrix embedding method \cite{wc}. The values of the constants involved in these simulations were, unless explicitly stated otherwise, $X_0 = 3$, $\alpha = 0.5$, $\beta =  2$, and $\sigma = 1.5$. We considered these sample paths on the unit interval, hence $T = 1$. {The number of replicates was 300 in all of the considered cases.} In what follows we present the dependencies of the estimators both on the true parameter {values ($H$, $\sigma$)} and on the {sample size ($n$)}. We have also checked for possible dependencies of the estimators of the Hurst index and {the diffusion coefficient} on the values of the other parameters of the considered equation, namely the drift coefficients $\alpha$ and $\beta$ and the initial condition $X_0$. No such dependencies of significant impact have been observed.

\subsection{Modeling of the Hurst index estimators}

\setcounter{figure}{0}

Figures 1 and 2 display, respectively, the dependence of the four estimators {of the Hurst index $H$  on its true value and on the sample size (length of the sample path) $n$. In Figure 1, the same sample sizes $n=2^{10}$ were used for all of the considered estimators, which does suggest that the estimators $\widehat H^{(4)}_n$ and $\widehat H^{(2)}_n$ would be a-priori less efficient. However, in practical applications the sample size is usually fixed, hence the motivation was to see what kind of performance the considered estimators would show given the exact same number of observations. In Figure 2, the value of the Hurst index was chosen as $H=0.75$. The values of $r_j$ were taken to be powers of 2 (more precisely, $r_j = 2^{j-1}$, $j=1,\,\ldots,\,l$) and, further, the values of $n_i$ were taken as $n_i = n / r_i$ where $n$ denotes the (fixed) maximum sample size length. The value of $l$ was (arbitrarily) taken to be $4$, as simulation results suggested that both considerably smaller (f.e., $2$) and considerably larger (f.e., $log_2{n}-1$) values yielded inferior performance. It does appear plausible that for much bigger sample sizes it might be beneficial to increase this value further, however in this study sample sizes exceeding 6400 points were not considered.} It can be seen that the performance of the estimator $\widehat H_n^{(4)}$ is slightly lacking compared to that of the other estimators{, which, despite imposing rather different requirements on the sample sizes, show similar precision.}

\begin{figure}[H]
\label{h_dep_h}
\centering
  \includegraphics[scale=0.6]{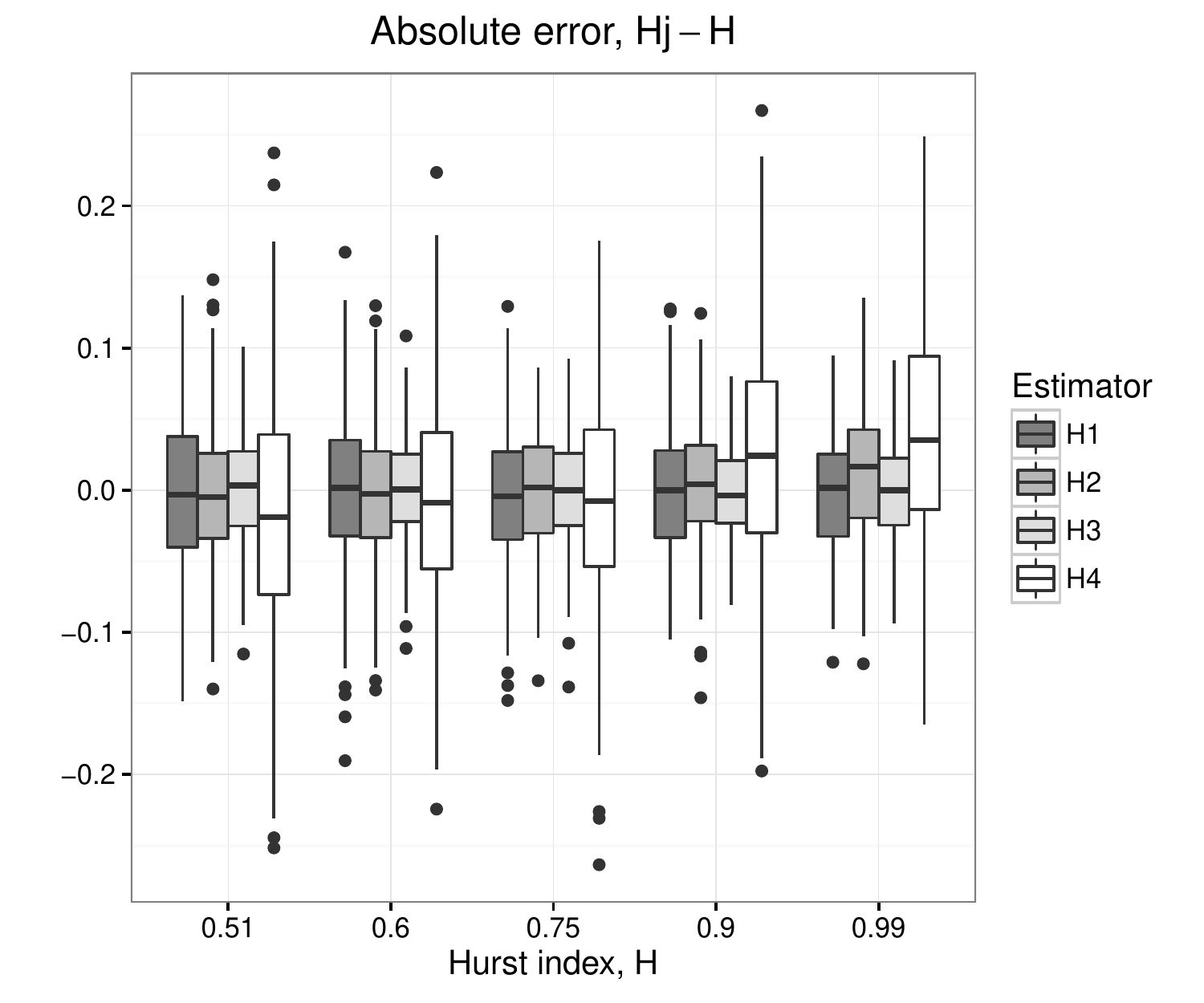}
  \caption{Dependence {of the absolute error on $H$}}
\end{figure}

\begin{figure}[H]
\label{h_dep_n}
\centering
  \includegraphics[scale=0.6]{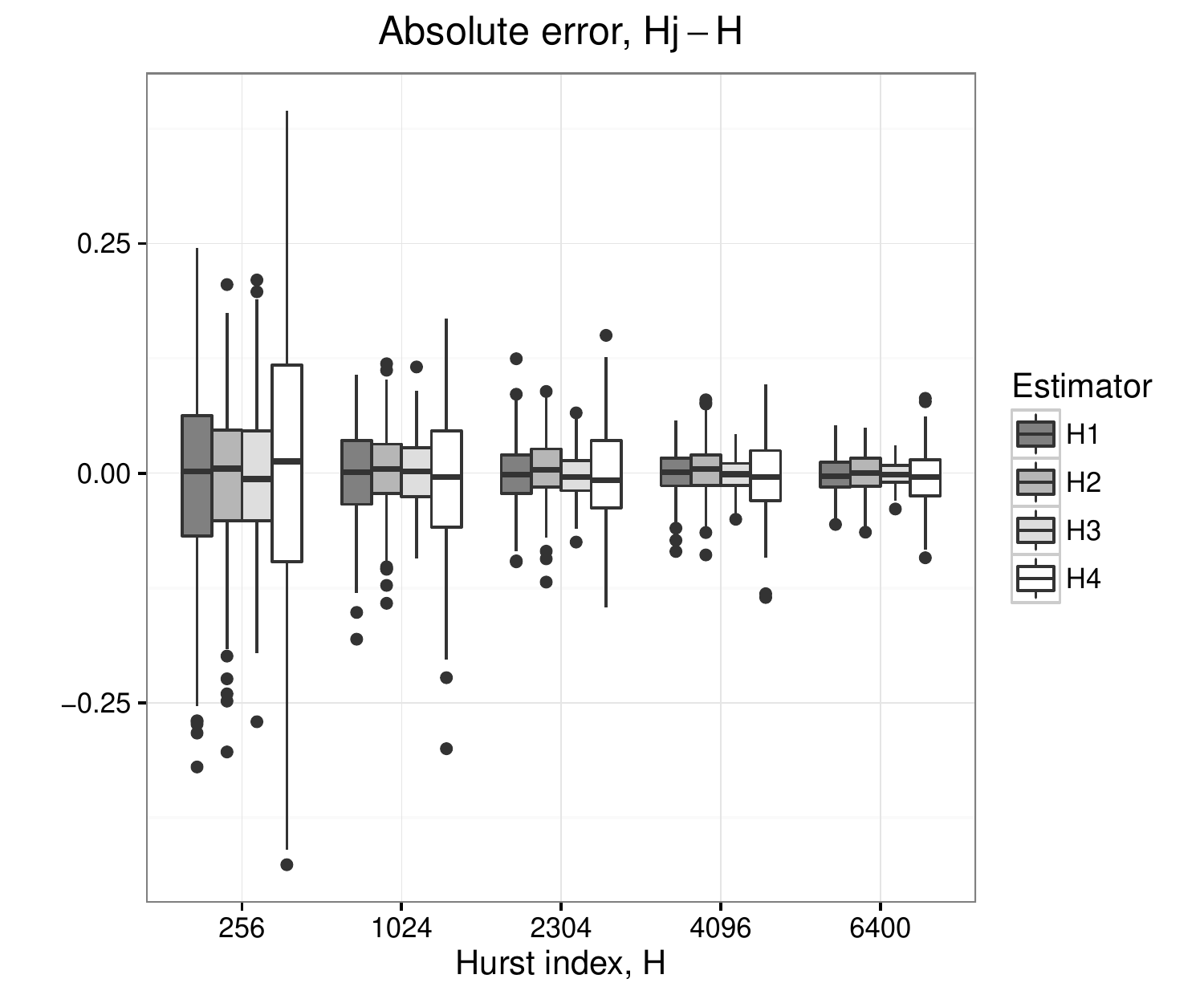}
  \caption{Dependence {of the absolute error on $n$}}
\end{figure}

\subsection{Modeling of {the diffusion coefficient} estimators}

In order to calculate the estimators $\widehat \s^{\,2}_{1,n}$, $\widehat \s^{\,2}_{2,n}$ and $\widehat \s^{\,2}_{3,n}$ we need to supply them with the estimated values of the Hurst index. In the Figures 3 and 4 presented below, {the diffusion coefficient} estimator $\widehat \s^{\,2}_{i,n}$, using the Hurst index estimator $\widehat H_n^{(j)}$, is denoted as `si\_hj', $i,\,j=1,\,2,\,3$. The estimator $\widehat \s^{\,2}_{4,n}$ is denoted as `s4'. The graphs present the relative differences, namely $(\widehat \s_{i,n} - \sigma) / \sigma$. {In Figure 3, the sample size was chosen to be $n=2^{10}$ for all of the considered estimators. In Figure 4, the value of the diffusion coefficient was chosen as $\sigma = 1$.}
It can be seen that the performance of all the considered estimators is roughly similar. The convergence rate of $\widehat \s^{\,2}_{4,n}$ appears slower, although it seems to perform better for the values of $\sigma$ close to zero. For the other estimators, it appears that using $\widehat H_n^{(3)}$ yields better numerical characteristics.

\begin{figure}[H]
\label{s_dep_s}
\centering
  \includegraphics[scale=0.6]{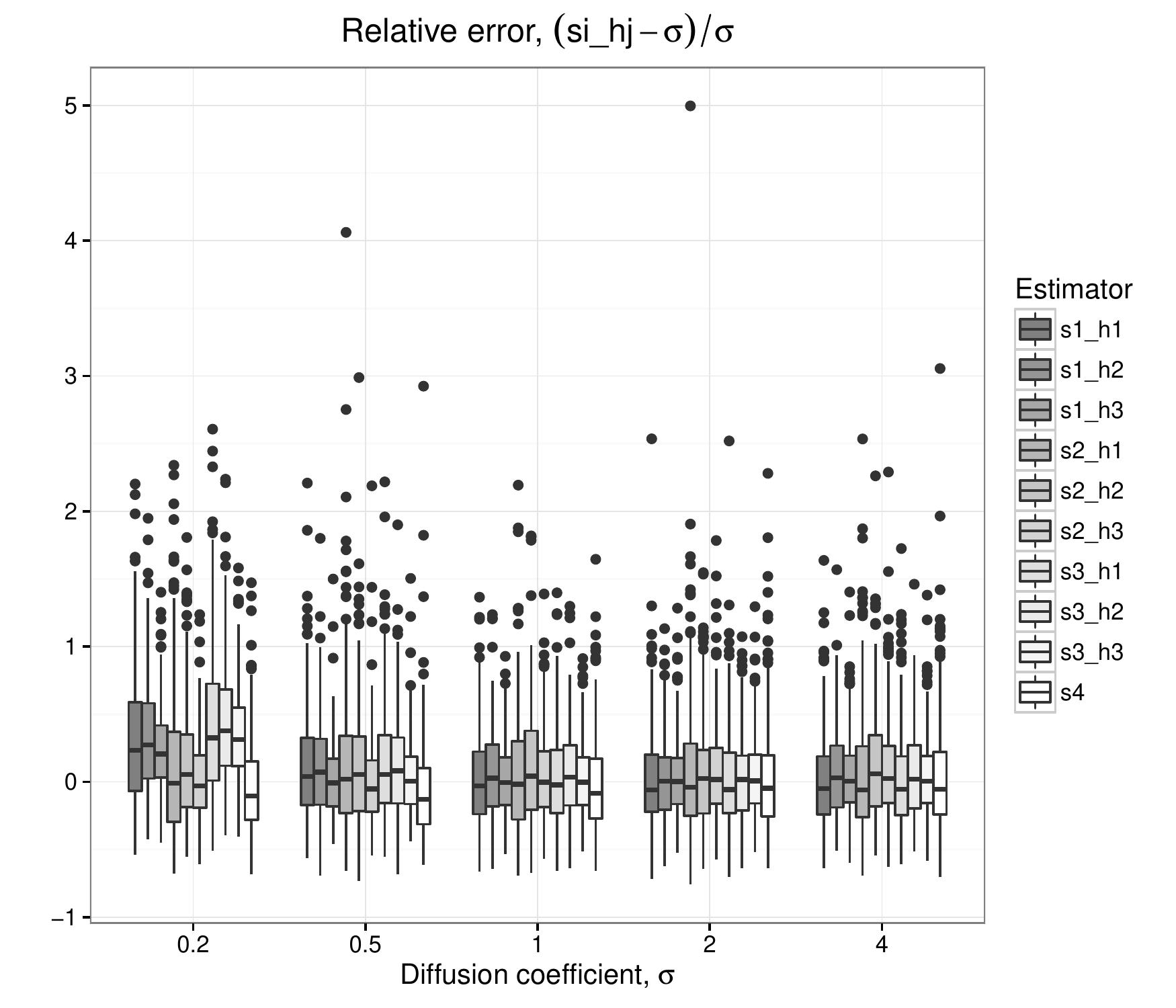}
  \caption{Dependence {of the relative error on $\sigma$}}
\end{figure}

\begin{figure}[H]
\label{s_dep_n}
\centering
  \includegraphics[scale=0.6]{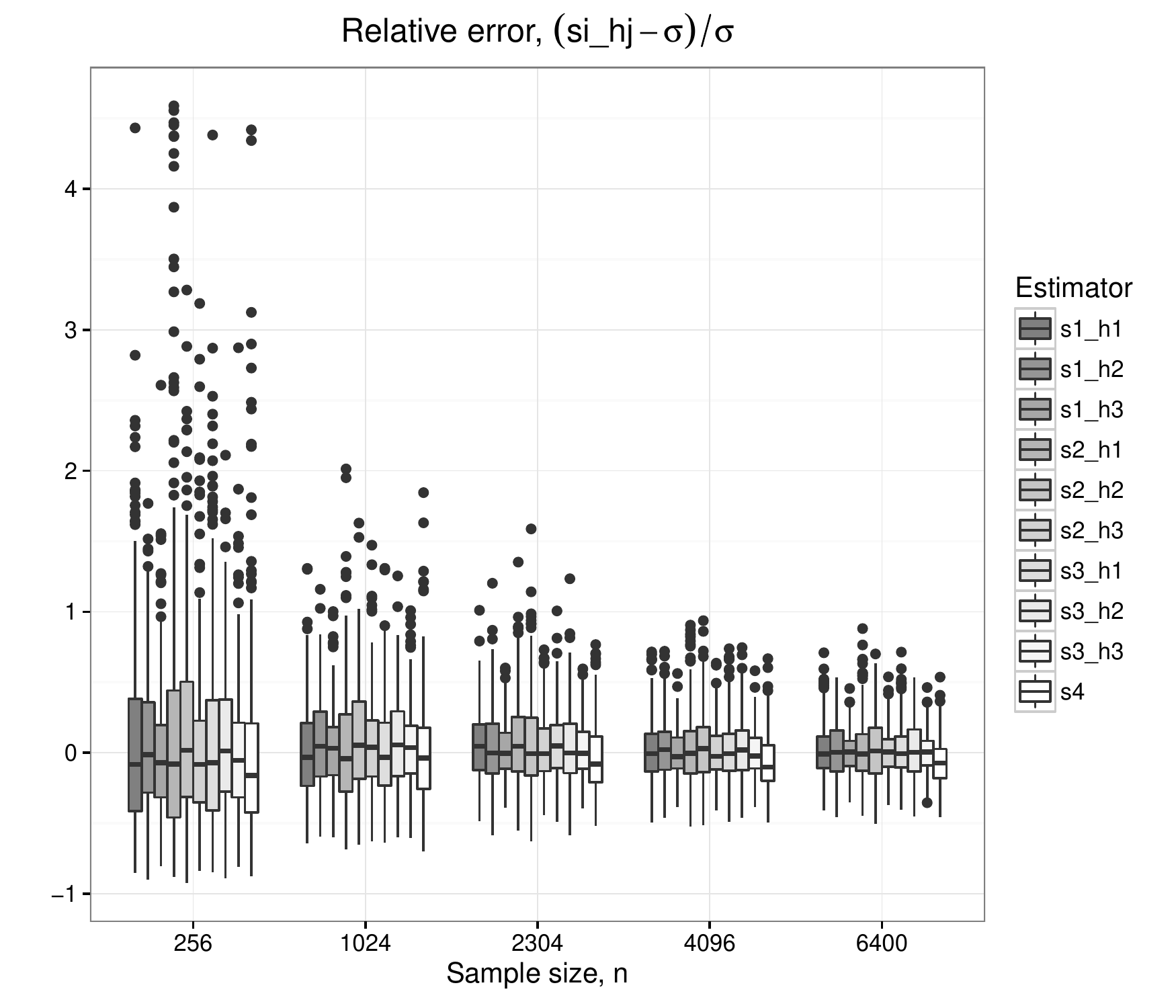}
  \caption{Dependence {of the relative error on $n$}}
\end{figure}

\section{Preliminaries}\label{s:preliminaries}

\subsection{Variation}\label{s:variation}

Let $p>0, -\infty<a<b<\infty$ be fixed and
$\mathcal{K}=\{\{x_0,\dots,x_n\}\mid a=x_0<x_1<\dots<x_n=b,n\ge 1\}$
denotes the set of all possible partitions of $[a,b]$. For any
$f:[a,b]\to\mathbb{R}$ define

\begin{gather*}
    v_{p}(f;[a,b])=\sup_{\varkappa\in\mathcal{K}}\sum_{k=1}^{n}\abs{f(x_k)-f(x_{k-1})}^{p},\qquad
    V_{p}(f;[a,b])=v_p^{\frac{1}{p}}(f;[a,b]).
\end{gather*}
If $v_p(f;[a,b])<\infty,$ $f$ is said to have a bounded
$p$-variation on $[a,b]$.\\

\noindent In the rest of the paper $\mathcal{W}_p([a,b])$ denotes the class of
functions on $[a,b]$ with bounded $p$--variation and
$C\mathcal{W}_p([a,b])=\{f\in\mathcal{W}_p([a,b])\mid f$ is
continuous$\}$. In case of a fixed interval $[a,b]$ we abbreviate the
notations and write $v_p(f)$, $V_p(f)$, etc. instead of
$v_p(f;[a,b])$, $V_p(f;[a,b])$.\\

\noindent Below we list several facts used in the sequel. For details we refer
the reader to \cite{DudleyNorvaisa-10}.

\begin{itemize}
  \item $p\ge1\Rightarrow f\mapsto V_p(f)$ is a semi-norm on
  $\mathcal{W}_p$.
  \item $f=const.\Longleftrightarrow V_p(f)=0$.
  \item $f\in \mathcal{W}_p\Rightarrow
  \sup_{x\in[a,b]}\abs{f(x)}<\infty$.
  \item $p\ge1,f\in \mathcal{W}_p\Rightarrow f\in \mathcal{W}_q$ for all $q\ge p$.
  \item $f,g\in \mathcal{W}_p\Rightarrow fg\in \mathcal{W}_p$.
  \item Let $f\in \mathcal{W}_q$ and $h\in \mathcal{W}_p$ with
        $p,q\in(0,\infty):1/p+\allowbreak 1/q>1.$ Then an integral
        $\int_a^b f\,\mathrm{d}h$ exists as the Riemann--Stieltjes  integral
        provided $f$ and $h$ have no common discontinuities. If the integral
        exists, the Love--Young inequality
        \begin{equation}\label{in:LoveYoung}
        \Bigg\vert \int\limits_a^bf\,\dr h-f(y)\big[ h(b)-h(a)
        \big]\Bigg\vert \ls C_{p,q} V_q\big(f\big)V_p\big(h\big)
        \end{equation}
        holds for all $y\in [a,b]$, where $C_{p,q}=\zeta(p^{-1}+q^{-1})$ and
        $\zeta(s)=\sum_{n\gs 1} n^{-s}$. Moreover,
        \[
        V_p\Bigg (\int\limits_a^{\bt} f \,\dr h\Bigg) \leq C_{p, q}
        V_{q,\infty}\big(f\big) V_p\big(h\big),
       \]
        where $V_{q,\infty}(f)=V_q(f)+\sup_{a\ls x\ls b}\vert
        f(x)\vert$. Note that $f\mapsto V_{q,\infty}(f)$ is a norm on ${\cal
        W}_q,q\gs 1$.
  \item $h\in C{\cal W}_p\Rightarrow g(y)=\int_a^y
f\,\dr h$, $y\in [a,b]$, is continuous.
\item {Let $\phi$ be a locally Lipschitz function and let $f\in\mathcal{W}_p([a,b])$. Then the composite {function} $\phi\circ f$ has bounded $p$-variation, that is, $\phi\circ f\in\mathcal{W}_p([a,b])$, where $(\phi\circ f)(x)=\phi(f(x))$.}
\end{itemize}

\noindent {The chain rule is based on the Riemann--Stieltjes integrals.}

\begin{thm}[Chain rule (see \cite{DudleyNorvaisa-10})]\label{t:salopek} Let $p\in[1;2)$ and $f = (f_1,\ldots, f_d): [a, b]\to \mathbb{R}^d$ be such a function that for each $k = 1,\ldots,d$, $f_k\in C \mathcal{W}_p([a,b])$.
Let $g: \mathbb{R}^d \to \mathbb{R}$ be a differentiable function
with locally Lipschitz partial derivatives $g^\prime_k$,
$k=1,\ldots,d$. Then each $g^\prime_l\circ f$ is Riemann-Stieltjes
integrable with respect to $f_k$ and
\[
(g\circ f)(b)-(g\circ f)(a)=\sum_{k=1}^d \int_a^b (g_{k}^\prime\circ
f)\, df_k.
\]
\end{thm}

\begin{prop}[Substitution rule (see \cite{DudleyNorvaisa-10})]\label{p: substitution rule} Let $f,g,h\in C \mathcal{W}_p([a,b])$, $1\ls p < 2$. Then
\[
\int_a^b f(x)\,d\bigg(\int_a^x g(y)\,dh(y)\bigg)=\int_a^b
f(x)g(x)\,dh(x).
\]
\end{prop}

\subsection{Several results on fBm}\label{s:fbm}

Recall that the fBm $(B^H_t)_{t\in[0,T]}$ with the Hurst index $H\in (0,1)$
is a real-valued continuous centered Gaussian process with the
covariance given by
\[
\E(B^H_t B^H_s)=\frac12\big(s^{2H}+t^{2H}-|t-s|^{2H}\big).
\]
In order to consider the strong consistency and asymptotic normality of
the given estimators we need several facts about $B^H$ (see  \cite{begyn2}, \cite{BLL}, \cite{Coeurjolly-01},  \cite{IL}).

\medskip\noindent\emph{\textbf{Limit results.}} For consideration of the asymptotic properties of the estimators $\widehat H^{(i)}_n$, $i=1,2$, we shall use the following results. Let
\[
\widehat V^{(i)B^H}_{n,T}=\frac{n^{2H-1}}{c_i} \sum_{k=1}^{n-1}\big(T^{-H}{\Delta^{(i)}
B^H(t^n_k)}\big)^2,\qquad H\neq \frac 12\,,\quad i=1,2.
\]
where $c_1=1,c_2=4-2^{2H}$. Then
\begin{gather*}
    \widehat V^{(i)B^H}_{n,T}\tends{n\to\infty} 1\qquad\text{ a.s.}\quad i=1,2; \\
    {\sqrt{n}\left(\widehat V^{(1)B^H}_{n,T}-1\right)\tendsd\mathcal{N}(0,\s^2_*(H)),\qquad H\in(0,3/4),}\\
    \sqrt{n}\left(
      \begin{array}{c}
        \widehat V^{(2)B^H}_{n,T}-1 \\
        \widehat V^{(2)B^H}_{2n,T}-1 \\
      \end{array}
    \right)\tendsd\mathcal{N}\left(\left(
                                     \begin{array}{c}
                                       0 \\
                                       0 \\
                                     \end{array}
                                   \right),
                                   \left(
               \begin{array}{cc}
                 \sigma^2(H) & \sigma_{1,2}(H) \\
                 \sigma_{1,2}(H) & \sigma^2(H)/2 \\
               \end{array}
             \right)
    \right)
\end{gather*}
with
\begin{gather*}
{\s^2_*(H)=2\bigg(1+2\sum_{j=1}^{\infty}\widehat\rho_H^2(j)\bigg), \qquad \widehat\rho_H(j)=-\frac12\big[-\abs{j-1}^{2H}+2\abs{j}^{2H}-\abs{j+1}^{2H}\big],}\\
    \sigma^2(H)=2\bigg(1+2\sum_{j=1}^{\infty}\rho_H^2(j)\bigg),\qquad \sigma_{1,2}(H)=\sum_{j\in\mathbb{Z}}\widetilde{\rho}_H^2(j)\,,\\
    \rho_H(j)=\frac{-6\abs{j}^{2H}-(\abs{j-2}^{2H}+\abs{j+2}^{2H})+4(\abs{j-1}^{2H}+\abs{j+1}^{2H})}{2(4-2^{2H})},\\
    \widetilde{\rho}_H(j)=\frac{
    -\abs{j-2}^{2H}+2\abs{j-1}^{2H}+\abs{j}^{2H}-4\abs{j+1}^{2H}+ \abs{j+2}^{2H}+2\abs{j+3}^{2H}-\abs{j+4}^{2H}
    }{2(4-2^{2H})2^{H}}\,.
\end{gather*}

\noindent In order to prove the asymptotic normality of the estimator $\widehat H^{(3)}_n$ we need the following result obtained in \cite{BLL}. Let $n_i=r_i n$, $i=1,\ldots,\ell$, where $r_i, n\in \mathbb{N}$, and $z_i$, $i=1,\ldots,\ell$, are defined in Theorem \ref{t:case1}. Then
\[
\frac 12 \sum_{i=1}^\ell \frac{z_i}{\sqrt{r_i}} \sqrt{n_i}\Big(\widehat V^{(i)B^H}_{n_i,T}-1\Big)\tendsd{} \mathcal{N}\bigg(0,\sigma_{2,\ell}^2\Big(\mathbf{r},\frac 12(\mathbf{z}/\sqrt{\mathbf{r}}\,)\Big)\bigg),
\]
where $\mathbf{r}=(r_1,\ldots,r_\ell)$, $\mathbf{z}=(z_1,\ldots,z_\ell)$,
\begin{align*}
\sigma_{2,\ell}^2(\mathbf{k},\mathbf{d})=&\sum_{i=1}^\ell \sum_{j=1}^\ell d_i d_j\rho_2(k_i,k_j),\qquad \mathbf{k}=(k_1,\ldots,k_\ell)\in\mathbb{N}^\ell,\quad \mathbf{d}=(d_1,\ldots,d_\ell)\in\mathbb{R}^\ell,\\
\rho_2(k_i,k_j)=&\frac{1}{\sqrt{k_i k_j}}\sum_{p=1}^{+\infty}c_{2p,2}^2\cdot (2p)! \bigg(\sum_{s=0}^{k_i-1}\sum_{r=-\infty}^{+\infty}\rho_{k_i,k_j}^{2p}(k_i r+k_j s)\bigg),\qquad c_{2p,2}=\frac{1}{(2p)!}\prod_{i=0}^{p-1}(2-2i),\\
\rho_{b,c}(x)=&\frac{1}{2(4-2^{2H})}\,(bc)^{-H}\big[-\abs{x}^{2H}+2\abs{x-b}^{2H} -\abs{x-2b}^{2H}+2\abs{x+c}^{2H}-4\abs{x+c-b}^{2H}\\
&+2\abs{x+c-2b}^{2H}-\abs{x+2c}^{2H}+2\abs{x+2c-b}^{2H}-\abs{x+2c-2b}^{2H}\big].
\end{align*}
If $k_i=k_j$ then
\[
\rho_2(k,k)=2\sum_{r=-\infty}^{+\infty}\rho_H^2(r) =2\bigg(1+2\sum_{j=1}^{\infty}\rho_H^2(j)\bigg).
\]

\medskip\noindent\emph{\textbf{Variation of $B^H$.}} It is known that almost all sample paths of $B^H$ are locally H\"older of order strictly less than $H$, $0<H<1$. To be more
precise, for all $0<\eps<H$ and $T>0$, there exists a nonnegative
random variable $G_{\eps,T}$ such that {$\mathbb{E}(\vert G_{\eps,T}\vert^p)<\infty$ for all $p\ge1$} and
\begin{equation}\label{in:fBmpokyt}
sup_{s,t\in[0;T]}\vert B^H_t -B^H_s\vert \ls G_{\eps,T}\vert
t-s\vert^{H-\eps}\qquad a.s.
\end{equation}
Thus $B^H\in C\mathcal{W}_{H_\eps}([0,T])$,
${H_\eps}=\frac{1}{H-\eps}$.

\medskip\noindent\emph{\textbf{The rate of convergence of the Hurst index.}}
\begin{thm}[\cite{km}, Theorem 3]\label{t:variacijos greitis}
For any $t\in[0;T]$ define
$r_{nt}=\left[\frac{tn}{T}\right],\rho_{nt}=\frac{r_{nt}}{n}T$ and
\[
{\widehat V_{nt}^{(2)B^H}=\frac{n^{2H-1}}{T^{2H-1}(4-2^{2H})} \sum_{k=i}^{r_{nt}}\big(\Delta^{(2)}B^H(t^n_k)\big)^2.}
\]
Then
\begin{equation}\label{e:e:variacijos_asimp}
  { \sup_{t\in[0;T]}\abs{\widehat V_{nt}^{(2)B^H}-\rho_{nt}}=O_\omega(n^{-1/2}\ln^{1/2}n),}
\end{equation}
where $O_\omega$ is defined in subsection \ref{s:gompertz}.
\end{thm}

\section{{Properties of the increments of the Gompertz diffusion process}}\label{s:gompertz}

{The fractional Gompertz diffusion process $X$ has the explicit solution given by
\[
X_t=\exp\bigg\{e^{-\b t}\ln x_0+\frac \a \b\big(1-e^{-\b t}\big)+\s\int_0^t e^{-\b(t-s)}dB^H_s\bigg\}, \qquad 0\ls t\ls T.
\]
Moreover, it is unique in $C{\cal W}_{\frac{1}{H-\eps}}([0,T])$ for all $\eps\in\left(0,H-\frac{1}{2}\right)$. The proof of this {can be found in the} Appendix.  Now we will {consider} the structure of increments of the Gompertz diffusion process.}\\

\noindent To avoid cumbersome expressions, we introduce the symbols
$O_{\omega}$ and $o_{\omega}$. Let $(Y_n)$ be a sequence of r.~v.s,
$\varsigma$ is an a.~s. non-negative r.~v. and
$(a_n)\subset(0,\infty)$ vanishes. $Y_n=O_\omega(a_n)$ means that
$\vert Y_n\vert\le \varsigma\cdot a_n$; $Y_n=o_\omega(a_n)$ means
that $\vert Y_n\vert\le \varsigma\cdot b_n$ with $b_n=o(a_n)$. In
particular, $Y_n=o_\omega(1)$ corresponds to the sequence $(Y_n)$
which tends to $0$ a.~s. as $n\to \infty$.

\begin{lem}\label{e:pokyciu_asimp} Suppose that $X$ satisfies $(\ref{e:diflygt})$, $\eps\in(0,H-\frac{1}{2})$ and partition $\pi_n$ of the interval $[0,T]$ is uniform. Then the following relations hold:
\begin{align}
\Delta X_{\tau^{m_n}_k}=&X_{\tau^{m_n}_{k-1}}\Big[\s\D B^H_{\tau^{m_n}_k} +O_\omega(d_n)\Big] =X_{\tau^{m_n}_{k-1}}O_\omega\big(d_n^{H-\eps}\big),\qquad k=1,\dots, m_n,\label{e:nario_asimp_0}\\
\Delta^{(2)} X_{\tau^{m_n}_k} =&X_{\tau^{m_n}_{k-1}}\big[\s\D^{(2)} B^H_{\tau^{m_n}_k} + O_\omega\big(d_n^{2(H-\eps)}\big)\big],\qquad k=2,\dots,
    m_n,\label{e:nario_asimp_1}
\end{align}
where $d_n=\tau^{m_n}_k-\tau^{m_n}_{k-1}$ and $d_n\to 0$ as $n\to\infty$. Moreover, $\E O_\omega(1)<\infty$.
\end{lem}

\noindent \proof For the sake of simplicity we will omit the index $m_n$ for the points
$\tau^{m_n}_k$. Let the sample path $t\mapsto X_t$ be
continuous. We first prove (\ref{e:nario_asimp_0}). Note that
\[
\Delta X_{\tau_k}=X_{\tau_k}-X_{\tau_{k-1}}\quad\mbox{and}\quad
X_{\tau_k}=X_{\tau_{k-1}}\exp\{\Delta Y_{\tau_k}\},
\]
where
\[
Y_t=e^{-\b t}\ln x_0 +\frac \a \b\big(1-e^{-\b t}\big)+\s\int_0^t e^{-\b(t-s)}dB^H_s.
\]
It is clear that
\begin{align*}
\Delta Y_{\tau_k}
=& e^{-\b \tau_{k-1}}\big(e^{-\b (\tau_k-\tau_{k-1})} -1\big)\bigg(\ln x_0 -\frac\a\b+\s \int_0^{\tau_{k-1}} e^{\b s} dB^H_s\bigg)\\
&+\s e^{-\b\tau_k}\int^{\tau_k}_{\tau_{k-1}} \big[e^{\b s}-e^{\b\tau_k}\big] dB^H_s +\s \int^{\tau_k}_{\tau_{k-1}} dB^H_s.
\end{align*}
From the Chain rule it follows that
\[
\int_0^t e^{\b s}dB^H_s=e^{\b t}B^H_t-\b \int_0^t e^{\b s}B^H_s ds.
\]
Thus
\[
\bigg\vert\int_0^{\tau_{k-1}} e^{\b s} dB^H_s\bigg\vert=\bigg\vert e^{\b \tau_{k-1}}B^H_{\tau_{k-1}} -\b\int_0^{\tau_{k-1}} e^{\b s}B^H_s\, ds\bigg\vert
\ls  e^{\vert\b\vert T}(\vert\b\vert T+1)\sup_{t\ls T}\big\vert B^H_t\big\vert\,.
\]
Provided
\[
{e^{-\b (\tau_k-\tau_{k-1})}=1+O(d_n),}
\]
it follows that
\[
Z_{k-1}:=e^{-\b \tau_{k-1}}\big(e^{-\b (\tau_k-\tau_{k-1})}-1\big)\bigg[\ln x_0 -\frac\a\b+\s\int_0^{\tau_{k-1}} e^{\b s} dB^H_s\bigg]
= O_\omega(d_n).
\]
Further
\begin{align*}
\bigg\vert\int_{\tau_{k-1}}^{\tau_k} \big[e^{\b s}-e^{\b \tau_k}\big] dB^H_s\bigg\vert\ls& C_{1,H_\eps} V_1\big(e^{\b \mathbf{\bigcdot}};[\tau_{k-1},\tau_k]\big) V_{H_\eps}\big(B^H;[\tau_{k-1},\tau_k]\big)\\
\ls& C_{1,H_\eps} e^{2\vert\b\vert T} \vert \b\vert (\tau_k-\tau_{k-1}) V_{H_\eps}\big(B^H;[\tau_{k-1},\tau_k]\big)\\
\ls& C_{1,H_\eps} G_{\eps,T} e^{2\vert\b\vert T} \vert \b\vert (\tau_k-\tau_{k-1})^{1+H-\eps} =O_\omega\big(d_n^{1+H-\eps}\big),
\end{align*}
since $\vert e^x-1\vert\ls \vert x\vert e^{\vert x\vert}$ for all $x\in \mathbb{R}$. Consequently,
\begin{align}\label{e:nario_asimp_2}
X_{\tau_k}=&X_{\tau_{k-1}}\exp\Big\{Z_{k-1} + O_\omega\big(d_n^{1+H-\eps}\big) +\s \Delta B^H_{\tau^n_k}\Big\}\nonumber\\
=&X_{\tau_{k-1}}\Big[1+Z_{k-1} + O_\omega\big(d_n^{1+H-\eps}\big) +\s \Delta B^H_{\tau^n_k}+O_\omega\big(d_n^{2(H-\eps)}\big)\Big]
\end{align}
and
\begin{equation}\label{e:pokyciu_asimp1}
{\Delta X_{\tau_k}=X_{\tau_{k-1}}\Big[Z_{k-1} +O_\omega\big(d_n^{2(H-\eps)}\big) +\s \Delta B^H_{\tau^n_k}\Big]=X_{\tau_{k-1}}O_\omega\big(d_n^{H-\eps}\big).}
\end{equation}
Since (see subsection \ref{s:fbm} and \cite{debicki2009estimates})
\[
\E\Big(\sup_{t\ls T}\big\vert B^H_t\big\vert\Big)^p<\infty\quad\mbox{and}\quad \E\vert G_{\eps,T}\vert^p<\infty
\]
for all $p\gs 1$, then $\E O_\omega(1)<\infty$.\\

\noindent Next we prove (\ref{e:nario_asimp_1}). Taking into account (\ref{e:nario_asimp_2}) and (\ref{e:pokyciu_asimp1}) we get
{\begin{align*}
\Delta^{(2)} X_{\tau_k}=&X_{\tau_k}\Big[Z_k +O_\omega\big(d_n^{2(H-\eps)}\big) +\s \Delta B^H_{\tau^n_{k+1}}\Big]
-X_{\tau_{k-1}}\Big[Z_{k-1} +O_\omega\big(d_n^{2(H-\eps)}\big) +\s \Delta B^H_{\tau^n_k}\Big]\\
=&X_{\tau_{k-1}}\Big[1+Z_{k-1} +O_\omega\big(d_n^{2(H-\eps)}\big) +\s \Delta B^H_{\tau^n_k}\Big]
\cdot\Big[Z_k +O_\omega\big(d_n^{2(H-\eps)}\big) +\s \Delta B^H_{\tau^n_{k+1}})\Big]\\
&-X_{\tau_{k-1}}\Big[Z_{k-1} +O_\omega\big(d_n^{2(H-\eps)}\big) +\s \Delta B^H_{\tau^n_k}\Big]\\
=&X_{\tau_{k-1}}\Big[(Z_k-Z_{k-1}) +O_\omega\big(d_n^{2(H-\eps)}\big) +\s \Delta^{(2)} B^H_{\tau^n_k}\Big] + X_{\tau_{k-1}} O_\omega\big(d_n^{2(H-\eps)}\big).
\end{align*}}
Since
\begin{align*}
Z_k-Z_{k-1}=&\big(e^{-\b \tau_k}-e^{-\b \tau_{k-1}}\big)\big(e^{-\b d_n}-1\big)\bigg[\ln x_0 -\frac\a\b+\s\int_0^{\tau_{k-1}} e^{\b s} dB^H_s\bigg]\\
&+\s e^{-\b \tau_k}\big(e^{-\b d_n}-1\big)\int^{\tau_k}_{\tau_{k-1}} e^{\b s} dB^H_s= O_\omega\big(d_n^{2(H-\eps)}\big),
\end{align*}
then
\begin{align*}
\bigg\vert\int_{\tau_{k-1}}^{\tau_k}  e^{\b s}\, dB^H_s\bigg\vert\ls& C_{1,H_\eps} V_{1,\infty}\big(e^{\b \mathbf{\bigcdot}};[\tau_{k-1},\tau_k]\big) V_{H_\eps}\big(B^H;[\tau_{k-1},\tau_k]\big)\\
\ls& 2C_{1,H_\eps} e^{\vert\b\vert T} V_{H_\eps}\big(B^H;[\tau_{k-1},\tau_k]\big)\\
\ls& 2C_{1,H_\eps} G_{\eps,T} e^{\vert\b\vert T} \vert \b\vert (\tau_k-\tau_{k-1})^{H-\eps} =O_\omega\big(d_n^{H-\eps}\big).
\end{align*}
Thus
\[
\Delta^{(2)} X_{\tau_k}=X_{\tau_{k-1}}\big[\s \Delta^{(2)} B^H_{\tau^n_k}+O_\omega\big(d_n^{2(H-\eps)}\big) \big].
\]

\section{Proofs of the main Theorems}

\subsection{Proof of Theorem \ref{t:case1}}

{1. The convergence of the statistics $\widehat H_n^{(1)}$ and $\widehat H_n^{(2)}$ considered in Theorem \ref{t:case1} follows from  Lemma \ref{e:pokyciu_asimp}. {Indeed}, the asymptotic{s} of {the} increments of the solution $X$ of the equation (\ref{e:diflygt}) {are} the same as {the} asymptotic{s} of {the} increments of the solution of the equation with polynomial drift in \cite{ksm}. Thus {in order to establish the convergence of the} estimator $\widehat H_n^{(1)}$ it suffices to repeat the proof of Theorem 2 in \cite{ksm}. Further, note that hypotheses $(H)$ and $(H_1)$ in \cite{ks} are satisfied for the solution of the equation (\ref{e:diflygt}), i.e.
\begin{align*}
\Delta X_{\tau^{m_n}_k}=&O_\omega\big(d_n^{H-\eps}\big),\qquad k=1,\dots, m_n,\\
\Delta^{(2)} X_{\tau^{m_n}_k} =&\s X_{\tau^{m_n}_{k-1}} \D^{(2)} B^H_{\tau^{m_n}_k} + O_\omega\big(d_n^{2(H-\eps)}\big),\qquad k=2,\dots, m_n.
\end{align*}
It follows from Lemma \ref{e:pokyciu_asimp} and {the} a.s. continuity of  $t\mapsto X_t$. Thus it suffices to apply Theorem 2.2 in \cite{ks}.}\\

\noindent 2. Now we prove the convergence of the statistic $\widehat H_n^{(3)}$.
The {proof presented below} follows the outline of the proof of Theorem 3.18 in \cite{BLL}. By Lemma \ref{e:pokyciu_asimp} we get
\begin{align}\label{e:kv_pokycio_asimp1}
\bigg(\frac{n^H}{\s T^H\sqrt{4-2^{2H}}}\bigg)^2 \frac{V^{(2)}_{n,T}}{n-1}
=&\bigg(\frac{n^H}{T^H\sqrt{4-2^{2H}}}\bigg)^2\frac{1}{n-1}\sum_{i=1}^{n-1} \bigg[ \big(\D^{(2)}_n B^H_i\big)^2 + O_\omega\big(n^{-3(H-\eps)}\big)\bigg]\nonumber\\
=&\frac{n}{n-1}\,\widehat V^{(2)B^H}_{n,T} +\frac{1}{4-2^{2H}}\,O_\omega\big(n^{-H+3\eps}\big)\tendsb 1.
\end{align}
Assume that $3\eps<H-1/2$. By (\ref{e:kv_pokycio_asimp1}) and Theorem \ref{t:variacijos greitis} we get
\begin{align*}
\ln\frac{V^{(2)}_{n,T}}{n-1}
=&-2 H\ln \frac nT+2\ln\big(\s \sqrt{4-2^{2H}}\big)+\ln\frac{n}{n-1}\\
&+\ln\bigg[ \big(\widehat V^{(2)B^H}_{n,T}-1\big)+1 +\frac{n-1}{n(4-2^{2H})}\,O_\omega\big(n^{-H+3\eps}\big)\bigg]\\
=&-2 H\ln \frac nT+2\ln\big(\s \sqrt{4-2^{2H}}\big)+\ln\frac{n}{n-1}\\
&+\ln\bigg[ O_\omega(n^{-1/2}\ln^{1/2}n)+1 +O_\omega\big(n^{-H+3\eps}\big)\bigg]\\
=&-2 H\ln \frac nT+2\ln\big(\s \sqrt{4-2^{2H}}\big) +O_\omega(n^{-1/2}\ln^{1/2}n).
\end{align*}
Thus
\[
\widehat H^{(3)}_n=-\frac 12\sum_{i=1}^\ell z_i\Big[-2H\ln \Big(\frac{n_i}{T}\Big)+2\ln\big(\s \sqrt{4-2^{2H}}\big)\Big] +O_\omega(n^{-1/2}\ln^{1/2}n).
\]
We will notice the following properties:
\[
\sum_{i=1}^\ell y_i=0,\qquad \sum_{i=1}^\ell z_i y_i=1,\qquad \sum_{i=1}^\ell z_i=\frac{\sum_{i=1}^\ell y_i}{\sum_{i=1}^\ell y_i^2}=0.
\]
Using those we get
\begin{align}\label{e:ivertinys}
\widehat H^{(3)}_n=& H\sum_{i=1}^\ell z_i \ln(r_in) -\ln\big(\s \sqrt{4-2^{2H}}\big)\sum_{i=1}^\ell z_i +O_\omega(n^{-1/2}\ln^{1/2}n)\nonumber\\
=&H\sum_{i=1}^\ell z_i \ln(r_in) +O_\omega(n^{-1/2}\ln^{1/2}n)\nonumber\\ =&H\sum_{i=1}^\ell z_i \big[y_i-y_i+\ln(r_in) \big] +O_\omega(n^{-1/2}\ln^{1/2}n)\nonumber\\
=&H+H \sum_{i=1}^\ell z_i\bigg[\ln n+\frac 1\ell \sum_{i=1}^\ell \ln r_i\bigg] +O_\omega(n^{-1/2}\ln^{1/2}n) =H+O_\omega(n^{-1/2}\ln^{1/2}n).
\end{align}
So the estimator $\widehat H^{(3)}_n$ is strongly consistent.\\

\noindent Now we prove the asymptotic normality of the estimator $\widehat H^{(3)}_n$. From (\ref{e:kv_pokycio_asimp1}) and (\ref{e:ivertinys}) it follows that
\begin{align*}
\widehat H^{(3)}_n=&H-\frac 12 \sum_{i=1}^\ell z_i\big(\widehat V^{(2)B^H}_{n_i,T}-1\big) +O_\omega\big(n^{-H+3\eps}\big).
\end{align*}
Thus
\begin{align*}
\sqrt{n}\big(\widehat H^{(3)}_n-H\big)=&-\frac 12 \sum_{i=1}^\ell \frac{z_i}{\sqrt{r_i}}\bigg[\frac{1}{\sqrt {r_in}} \sum_{k=1}^{r_in-1}\bigg(\bigg(\frac{(r_in)^H}{T^H\sqrt{4-2^{2H}}}\,\Delta^{(2)}
B^H(t^{r_in}_k) \bigg)^2-1\bigg)\bigg]\\
&+O_\omega\big(n^{1/2-H+3\eps}\big)
\end{align*}
and we obtain the asymptotic normality of the estimator $\widehat H^{(3)}_n$ by the application of the limit results from subsection \ref{s:fbm}.\\

\noindent 3. It remains to determine the convergence of $\widehat H_n^{(4)}$. Denote
\[
R^{2,n}(X) = \frac{1}{n^4-2}\sum_{k=1}^{n^4-2}\frac{\vert \D^{(2)}X(\tau^{m_n}_k) + \D^{(2)}X(\tau^{m_n}_{k+1})\vert}{\vert \D^{(2)}X(\tau^{m_n}_k)\vert + \vert \D^{(2)}X(\tau^{m_n}_{k+1}) \vert}\,,\qquad \L_2(H)=\E \frac{\vert \D^{(2)} B^H_1 + \D^{(2)} B^H_2 \vert}{\vert \D^{(2)} B^H_1 \vert + \vert \D^{(2)} B^H_2 \vert}\,,
\]
where $\D^{(2)}B^H_j=B^H(j+1)-2B^H(j)+B^H(j-1)$, $j=1,2$. This statistic was introduced in \cite{bs}. Further on, we will require the following lemma which is a simple modification of Lemma 3.1 in \cite{bs}. In this lemma we have lifted the requirement for the random variables $Z_1$ and $Z_2$ to be independent. This became possible due to the application of less precise estimators of the partial derivatives.

\begin{lem}\label{lem1} Let $\psi(x_1,x_2)=\frac{\vert x_1+x_2\vert}{\vert x_1\vert+\vert x_2\vert}$, $x_1,x_2\in\mathbb{R}$, and let $(Z_1,Z_2)$ be a Gaussian vector with zero mean and variance $\E Z_i^2=1$,  $i=1,2$. Then for any r.~v. $\xi_i$, $i=1,2$ with finite second moments we have
\begin{equation}\label{e:atstumas}
\E\big\vert\psi(Z_1+\xi_1,Z_2+\xi_2)-\psi(Z_1,Z_2)\big\vert\ls 23\max_{i=1,2}\root 3\of{\E\xi_i^2}.
\end{equation}
\end{lem}
Let us proceed to the following claim.

\begin{prop}\label{prop1}
Let  $X$ be the solution of the fractional Gompertz SDE observed at times $\tau^{m_n}_k=\frac{k}{n^4}T$, $k=0,1,\ldots,n^4$. Then
\[
R^{2,n}(X) \tendsb \L_2(H)\quad\mbox{as } n\to \infty\quad\mbox{for } H\in(1/2,1)\,.
\]
\end{prop}
\proof For the sake of simplicity we will omit the index $m_n$ for the points $\tau^{m_n}_k$ and denote $d_n=\frac{T}{n^4}$. From Lemma \ref{e:pokyciu_asimp} it follows that
\[
\Delta^{(2)} X_{\tau_k} + \Delta^{(2)} X_{\tau_{k+1}}=\s X_{\tau_k} \big[\D^{(2)} B^H_{\tau_k} + \D^{(2)} B^H_{\tau_{k+1}}+ \zeta_1+\zeta_2\big]
\]
for every $\eps\in(0,H-\frac{1}{2})$, where
\[
\zeta_1=O_\omega\big(d_n^{2(H-\eps)}\big),\qquad \zeta_2=O_\omega\big(d_n^{H-\eps}\big)\big[\s\D^{(2)} B^H_{\tau_{k+1}} + O_\omega\big(d_n^{2(H-\eps)}\big)\big]=O_\omega\big(d_n^{2(H-\eps)}\big).
\]
Therefore
\[
R^{2,n}(X)
=\frac{1}{n^4-2}\sum_{k=1}^{n^4-2}\frac{\vert Z_1+Z_2+\xi_1+\xi_2 \vert}{\vert Z_1+\xi_1 \vert + \vert  Z_2+\xi_2 \vert}\,,
\]
where
\begin{align*}
Z_1=&\frac{1}{d_n^H\sqrt{4-2^{2H}}}\,\D^{(2)} B^H_{\tau_k}\,,\qquad Z_2=\frac{1}{d_n^H\sqrt{4-2^{2H}}}\,\D^{(2)} B^H_{\tau_{k+1}}\,,\\ \xi_1=&\frac{\zeta_1}{d_n^H\sqrt{4-2^{2H}}}\,,\qquad \xi_2=\frac{\zeta_2}{d_n^H\sqrt{4-2^{2H}}}
\end{align*}
and
\[
\E Z^2_1=\frac{n^{8H}}{T^{2H}(4-2^{2H})}\,\E\big(\D^{(2)} B^H_{\tau_k}\big)^2 =1.
\]
Let us apply Lemma \ref{lem1}. From the inequality (\ref{e:atstumas}) it follows that
\begin{align*}
\E\big\vert R^{2,n}(X)-R^{2,n}(B^H)\big\vert= \bigg(\frac{d_n^{2(H-2\eps)}}{4-2^{2H}}\bigg)^{1/3}\root 3\of{\E O_\omega(1)} =d_n^{2(H-2\eps)/3}\root 3\of{\E O_\omega(1)}.
\end{align*}
Then the Chebyshev's inequality yields
\begin{align*}
\pr\big(\big\vert  R^{2,n}(X)-R^{2,n}(B^H)\big\vert > n^{-\b}\big)\ls& n^\b d_n^{2(H-2\eps)/3}\root 3\of{\E O_\omega(1)}\\
<&T^{2(H-2\eps)/3}n^{{\b-8(H-2\eps)/3}}\root 3\of{\E O_\omega(1)}
\end{align*}
for $\eps\in(0,(H-\frac{1}{2})/2)$, $0<\b<1/3$ and
\begin{align*}
\sum_{n=1}^{\infty}\pr\big(\big\vert  R^{2,n}(X)-R^{2,n}(B^H)\big\vert > n^{-\b}\big)\ls \root 3\of{\E O_\omega(1)}\sum_{n=1}^{\infty}n^{\b-8(H-2\eps)/3}<\infty.
\end{align*}
According to the Borel--Cantelli Lemma,
\[
\pr\Big(\limsup_{n\to\infty}\big\{\big\vert  R^{2,n}(X)-R^{2,n}(B^H)\big\vert > n^{-\b}\big\}\Big)=0
\]
which implies that $R^{2,n}(X)\tendsb R^{2,n}(B^H)$, $n\to\infty$.\\

\noindent The convergence $R^{2,n}(B^H)\tendsb{} \Lambda_2(H)$, $n\to\infty$ is established in \cite{bs}
and holds for  $H\in(0;1)$. Clearly, provided $R^{2,n}(X)\tendsb R^{2,n}(B^H)$ and $R^{2,n}(B^H)\tendsb \Lambda_2(H)$, $n\to\infty$, it follows that $R^{2,n}(X)\tendsb\Lambda_2(H)$, $n\to\infty$, which completes the proof.\qquad$\square$
\\

\noindent The estimator $\widehat H^{(4)}_n$ based on $R^{2,n}(X)$ can be obtained using the approximation formula provided in Remark 4.3 \cite{bs}.

\subsection{Proof of Theorem \ref{t:case2}}

The proof of {the} convergence of $\widehat{\s}^{\,2}_{2,n}$ is analogous to that of $\widehat{c}^{\,2}_n$ in \cite{ksm}. Let us prove that $\widehat{\s}^{\,2}_{1,n}\tendsb \s^2$, as $n\to\infty$.
Suppose that $d_n=\frac Tn$. From Lemma \ref{e:pokyciu_asimp} it follows that
\begin{align*}
d_n^{-2H}n^{-1} V^{(1)}_{n,T}
=&\s^2d_n^{-2H}n^{-1}\sum_{i=1}^n \big(\D B^H_{t^n_k}\big)^2 +d_n^{-2H}O_\omega\big(d_n^{1+H-\eps)}\big)\nonumber\\
=&\s^2 \widehat V^{(1)B^H}_{n,T}+O_\omega\big(d_n^{1-H-\eps}\big).
\end{align*}
Since
\begin{equation}\label{e:asimpt}
\widehat V^{(1)B^H}_{n,T} \tendsb 1\quad\mbox{and}\quad \frac{n^{2(H-\widehat H_n)}}{T^{2(H-\widehat H_n)}} =\exp\left\{o_\omega(\phi(n))\ln\left(\frac{n}{T}\right)^2\right\}\tend{} 1,
\end{equation}
it can be concluded that
\[
\widehat{\s}^{\,2}_{1,n}=\frac{n^{2 \widehat H_n-1}}{T^{2 \widehat
H_n}}\,V^{(1)}_{n,T}\tendsb{} \s^2\,.
\]

\noindent Further, let us prove that $\widehat{\s}^{\,2}_{3,n}\tendsb \s^2$. Denote $d_n=\tfrac Tn$.  By (\ref{e:asimpt}) it suffices to show that
\[
\widetilde{\s}^{\,2}_n=\frac{\sum_{k=1}^n\big(\D^{(1)} X_{t^n_k}
\big)^2}{d_n^{2H}  \sum_{k=1}^n X^2_{t^n_{k-1}}}\tendsb{} \s^2.
\]
Notice that
\begin{align*}
\frac{\sum_{k=1}^n\big(\D^{(1)} X_{t^n_k}
\big)^2}{d_n^{2H}  \sum_{k=1}^n X^2_{t^n_{k-1}}}=\frac{d_n^{1-2H}\sum_{k=1}^n\big(\D^{(1)} X_{t^n_k}
\big)^2}{d_n \sum_{k=1}^n X^2_{t^n_{k-1}}}
\end{align*}
and
\begin{equation}\label{e:kv_pokycio_asimp}
d_n^{1-2H}\sum_{k=1}^n\big(\D^{(1)} X_{t^n_k}
\big)^2 =\s^2d_n^{1-2H}\sum_{k=1}^n X^2_{t^n_{k-1}}\big(\D B^H_{t^n_k}\big)^2 +O_\omega\big(d_n^{1-H-\eps}\big).
\end{equation}
In order to estimate (\ref{e:kv_pokycio_asimp}), observe that
\begin{eqnarray*}
d_n^{1-2H}\sum_{k=1}^n X^2_{t^n_{k-1}}\big(\D B^H_{t^n_k}\big)^2 =\int_0^T X^2_t\,d \widehat V^{(1)B^H}_{nt}
\end{eqnarray*}
and (see Theorem 7 in \cite{km})
\begin{align*}
{d_n^{1-2H}\sum_{k=1}^n X^2_{t^n_{k-1}}\big[\big(\D B^H_{t^n_k}\big)^2-\E\big(\D B^H_{t^n_k}\big)^2\big]}
= \int_0^T X^2_t\,d\big(V^{(1)B^H}_{nt}-\E V^{(1)B^H}_{nt}\big)\tendsb 0.
\end{align*}
Since
\[
d_n^{1-2H}\sum_{k=1}^n X^2_{t^n_{k-1}}\E \big(\D B^H_{t^n_k}\big)^2=d_n \sum_{k=1}^n X^2_{t^n_{k-1}}\tendsb \int_0^T X^2_t\,dt,
\]
then
\[
\widetilde \s^2_n\tendsb \s^2.
\]

\subsection{The convergence rate of $\widehat H_n^{(i)},$ $i=1,2,3.$}

Theorem \ref{t:case2} makes use of the conditions $\widehat H_n=H+o_\omega(\phi(n))$,  $\phi(n)=o\left(\frac{1}{\ln n}\right)$ for strong consistency. Let us show that this indeed holds {for $\widehat H_n^{(i)}$, $i=1,2,3$.}\\

\noindent\emph{The convergence rate of $\widehat H_n^{(1)}$}. From Lemma \ref{e:pokyciu_asimp} and the proof of Theorem 2 in \cite{ksm} it follows that
\[
\widehat H_n^{(1)}=\widetilde H_n+O_\omega\big(n^{-H+3\eps}\big),
\]
where
\begin{equation}\label{e:greitis}
\widetilde H_n=\frac{1}{2}-\frac{1}{2\ln
    2}{\ln\left(\frac{\widehat V^{(2)B^H}_{2n,T}}{2^{2H-1}\widehat V^{(2)B^H}_{n,T}}\right)} =H-\frac{1}{2\ln
    2}{\ln\left(\frac{\widehat V^{(2)B^H}_{2n,T}}{\widehat V^{(2)B^H}_{n,T}}\right).}
\end{equation}
It suffices to consider the convergence rate of the logarithmic term in the equation (\ref{e:greitis}). {Using} Theorem \ref{t:variacijos greitis} we get
\[
{\ln\left(\frac{\widehat V^{(2)B^H}_{2n,T}}{\widehat V^{(2)B^H}_{n,T}}\right)} =\ln\bigg(\frac{1+O_\omega((2n)^{-1/2}\ln^{1/2} (2n))}{1+O_\omega(n^{-1/2}\ln^{1/2} n)}\bigg)=\ln\big(1+o_\omega(n^{-1/2}\ln n)\big)=o_\omega\big(n^{-1/2}\ln n\big).
\]
Then the statistic $\widetilde H_n$ has the convergence rate of $o_\omega(n^{-1/2}\ln n)$. Consequently, $\widehat H_n^{(1)}$ satisfies the required condition  if $\eps<(H-1/2)/3$.\\

\noindent\emph{The convergence rate of $\widehat H_n^{(2)}$}. Denote
\[
S_{n,T}:=\frac{2}{nk_n^{2H-1}} \sum_{k=2}^n\frac{\big(\D^{(2)} X_{t^n_k}\big)^2}{W_{n,k-1}}\,.
\]
Then
\[
\widehat H^{(2)}_n=H+\frac{\ln S_{n,T}}{2\ln k_n}=H+\frac{\ln S_{n,T}}{4\ln n}\,.
\]
Proceeding along the lines of the proof of Theorem 2.2 from \cite{ks}, it can be concluded that
\begin{align*}
S_{n,T}=&\frac{\widetilde V^{B^H}_{n,T}+O_\omega\big(n^{-(H-3\eps)}\big)}{1 +O_\omega\left(\sqrt{\frac{\ln n}{k_n}}\,\right)+O_\omega\left(\frac{m_n^{2\eps}}{n^{H-\eps}}\right)} =\frac{\widetilde V^{B^H}_{n,T}+O_\omega\big(n^{-(H-3\eps)}\big)}{1 +O_\omega\left(\sqrt{\frac{\ln n}{n^2}}\,\right)+O_\omega\left(\frac{1}{n^{H-7\eps}}\right)}\,.
\end{align*}
If $\eps<(H-1/2)/7$, then
\[
S_{n,T}=\frac{1+O_\omega\Big(\sqrt{\frac{\ln n}{n}}\,\Big)}{1 +O_\omega\left(\frac{1}{n^{H-7\eps}}\right)}=1+{O_\omega\bigg(\frac{\ln n}{\sqrt{n}}\,\bigg).}
\]
Hence $\widehat H_n^{(2)}=H+o_\omega(\frac{1}{\ln n})$, if $\eps<(H-1/2)/7$.\\

\noindent\emph{The convergence rate of $\widehat H_n^{(3)}$} was obtained in the proof of Theorem \ref{t:case1}.

{\section*{Appendix}\label{s:solution}}

{\subsubsection*{Auxiliary results}}

\noindent {{Firstly,} we consider a non-random integral equation
\begin{equation}\label{e:diflygt1}
x_t=x_0+\int_0^t(\alpha x_s-\beta x_s\ln x_s)\, ds + \s \int_0^t x_s\, d h_s,\qquad x_0\gs0,\quad \b\neq 0,\quad 0\ls t\ls T,
\end{equation}
where $h\in C\mathcal{W}_p([0,T])$, $1<p<2$, and prove two auxiliary theorems used in the sequel.}\\

\begin{thm} The function
\begin{equation}\label{e:solution}
x_t= \exp\bigg\{e^{-\b t}\ln x_0+\frac \a \b\big(1-e^{-\b t}\big)+\s\int_0^t e^{-\b(t-s)}dh_s\bigg\},\qquad t\in[0,T],
\end{equation}
is an element of $C\mathcal{W}_p([0,T])$, $1< p<2$, and  satisfies {the} equation (\ref{e:diflygt1}).
\end{thm}
\proof We show that $x\in C\mathcal{W}_p([0,T])$, $1< p<2$.
Let
\[
z_t=e^{-\b t}\ln x_0+\frac \a \b\big(1-e^{-\b t}\big)+\s\int_0^t e^{-\b(t-s)}dh_s.
\]
It is evident that $z\in C\mathcal{W}_p([0,T])$, $1< p<2$. Thus by {the} property of composition of functions (see subsection \ref{s:variation}) we get $x\in C\mathcal{W}_p([0,T])$, $1< p<2$.\\

\noindent Now we verify that the function (\ref{e:solution}) satisfies
(\ref{e:diflygt1}). This statement can be checked by the application of the Chain rule and the Substitution rule. Namely, let $F(t,x,y)=\exp\{e^{-\b t}(\ln x_0+\a x+\s y)\}$ and denote
\[
A_t=\int_0^t e^{\b s} ds,\qquad C_t=\int_0^t e^{\b s} dh_s.
\]
{Note  that $x_t = F(t,A_t,C_t)$ and
\begin{align}\label{e:chain_rule}
F(t,A_t,C_t) =& F(0; 0; 0) +\int_0^t \partial_t F(s,A_s,C_s)dA_s +\int_0^t \partial_x F(s,A_s,C_s)dC_s\nonumber\\
&+\int_0^t \partial_y F(s,A_s,C_s)ds.
\end{align}
It follows from (\ref{e:chain_rule}) and Proposition \ref{p: substitution rule}
\begin{align*}
x_t=&x_0-\b \int_0^t x_s \ln x_s\, ds+\a \int_0^t x_s e^{-\b s} dA_s +\s \int_0^t x_s e^{-\b s} dC_s\\
=&x_0-\b \int_0^t x_s \ln x_s\, ds+\a \int_0^t x_s\,ds +\s \int_0^t x_s \,dh_s,
\end{align*}
since}
\[
dA_s=e^{\b s} ds,\qquad dC_s=e^{\b s} dh_s.
\]

{\begin{thm} The integral equation (\ref{e:diflygt1}) has a unique solution in $C\mathcal{W}_p([0,T])$, $1< p<2$.
\end{thm}
\proof We have already shown that at least one solution $x\in
C\mathcal{W}_p([0,T])$ exists. Assume it is not unique and $y\in
C\mathcal{W}_p([0,T])$ is a different one.}\\

\noindent Further, one can find a set of points
$0=\tau_0<\tau_1<\tau_2<\cdots<\tau_n=T$ which satisfies
\[
V_p(h;[\tau_{k-1},\tau_k])\ls \big(4 \vert \s\vert C_{p,p} \big)^{-1}
\]
for all $k$. Assume we have proved that
$x_{\tau_{k-1}}=y_{\tau_{k-1}}$.\\

\noindent {Using the} well known inequality $\ln(1+x)\ls  x$, $x>0$, we get
\begin{align*}
\vert \ln x_s-\ln y_s\vert=\bigg\vert \ln\left(1+\frac{y_t-x_t}{x_t}\right)\bigg\vert \ls \bigg\vert \frac{y_t-x_t}{x_t}\bigg\vert \ls L_{x,T}\abs{y_t-x_t}
\end{align*}
and
\[
\vert \ln x_s\vert \ls \Big\vert \ln \Big(\max_{0\ls t\ls T}x_s\Big)\Big\vert=:\widehat L_{x,T}
\]
where $L_{x,T}=(\min_{0\ls t\ls T}\abs{x_t})^{-1}>0$. Then
\begin{align*}
V_{p,\infty}(x-y;[\tau_{k-1},\tau_k]) =&V_{p,\infty}(x-y-(x_{\tau_{k-1}}-y_{\tau_{k-1}});[\tau_{k-1},\tau_k])\\
\ls & 2\vert \a\vert\int_{\tau_{k-1}}^{\tau_k} \vert
x_t-y_t\vert\, dt+2\vert \b\vert \int_{\tau_{k-1}}^{\tau_k} \vert
x_t\ln x_t-y_t\ln y_t\vert\,dt\\
&+2\vert \s\vert C_{p,p} V_{p,\infty}(x-y;[\tau_{k-1},\tau_k]) V_p(h;[\tau_{k-1},\tau_k])\\
\ls&2\big(\vert \a\vert+\vert\b\vert \widehat L_{x,T}+ \vert
\b\vert L_{x,T}\big) \int_{\tau_{k-1}}^{\tau_k} \vert x_t-y_t\vert\, dt\\
&+2\vert \s\vert C_{p,p} V_{p,\infty}(x-y;[\tau_{k-1},\tau_k])
V_p(h;[\tau_{k-1},\tau_k])
\end{align*}
and
\begin{align*}
V_{p,\infty}(x-y;[\tau_{k-1},\tau_k])\ls& 4 \big(\vert \a\vert+\vert
\b\vert \widehat L_{x,T}+ \vert
\b\vert L_{x,T}\big) \int_{\tau_{k-1}}^{\tau_k} \vert x_t-y_t\vert\, dt \\
\ls& 4 \big(\vert \a\vert+\vert\b\vert \widehat L_{x,T}+ \vert\b\vert L_{x,T}\big) \int_{\tau_{k-1}}^{\tau_k}
V_{p,\infty}(x-y;[\tau_{k-1},t])\, dt.
\end{align*}
Therefore by Gronwall's inequality
$V_{p,\infty}(x-y;[\tau_{k-1},\tau_k])=0$ and we can conclude that
$x=y$ on $[\tau_{k-1},\tau_k]$. Since $x_{\tau_0}=x_0=y_{\tau_0}$
the claim of the theorem follows from the repetitive application of the reasoning {explained above}.

\subsubsection*{The solution of SDE}

Since almost all sample paths of $B^H$, $1/2<H<1$, are continuous
and have bounded $H_\eps=\frac{1}{H-\eps}$-variation,
$\eps\in(0,H-1/2)$, the pathwise Riemann-Stieltjes  integral
$\int_0^t X_s\, d B^H_s$ exists for $X\in
C\mathcal{W}_{H_\eps}([0,T])$. So SDE (\ref{e:diflygt}) is well
defined for almost all $\omega$ and the obtained result for a non-random
integral equation can be applied to an equation driven by fBm.

\begin{thm} Suppose that $X_0>0$ and $m\gs 2$. The stochastic process
\[
X_t=\exp\bigg\{e^{-\b t}\ln x_0+\frac \a \b\big(1-e^{-\b t}\big)+\s\int_0^t e^{-\b(t-s)}dB^H_s\bigg\}, \qquad \b\neq 0,\quad 0\ls t\ls T.
\]
for almost all $\omega$ belongs to  $C\mathcal{W}_{H_\eps}([0,T])$
and is the unique solution of (\ref{e:diflygt}).
\end{thm}

\noindent {\textbf{Acknowledgment.} The authors would like to thank the referees for many valuable comments {which allowed us to} improve this paper.}

\end{document}